\documentclass[12pt,a4]{amsart}
\usepackage{amsmath,amssymb,amsthm,comment}
\usepackage{graphicx}
\theoremstyle{plain}
\newtheorem{thm}{Theorem}[section]
\newtheorem{lem}[thm]{Lemma}
\newtheorem{cor}[thm]{Corollary}

\newtheorem{defn}[thm]{Definition}

\newtheorem{example}[thm]{Example}

\newtheorem{rem}[thm]{Remark}

\newcommand{\R}{\mathbb{R}}

\newcounter{mycounter}

\usepackage[dvipdfm,
bookmarks=true,bookmarksnumbered=true,
bookmarksopen=true]{hyperref}
\begin{document}

\title{There are no 76 equiangular lines in $\mathbb{R}^{19}$}
\author{Wei-Hsuan Yu}
\address{Department of Mathematics, Michigan State University,
619 Red Cedar Road, East Lansing, MI 48824}
\email{u690604@gmail.com}
\date{}
\maketitle

\begin{abstract}
Maximum size of equiangular lines in $\mathbb{R}^{19}$ has been known in the range between 72 to 76 since 1973. Acoording to the nonexistence of strongly regular graph $(75,32,10,16)$ \cite{aza15}, Larmen-Rogers-Seidel Theorem \cite{lar77} and Lemmen-Seidel bounds on equiangular lines with common angle $\frac 1 3$ \cite{lem73}, we can prove that there are no 76 equiangular lines in $\mathbb{R}^{19}$. As a corollary, there is no strongly regular graph $(76,35,18,14)$. Similar discussion can prove that there are no 96 equiangular lines in $\mathbb{R}^{20}$. 
\end{abstract}

\section{Introduction}

A set of lines in $\R^n$ is called equiangular if the angle between each
pair of lines is the same. We are interested in the upper bounds on the number of
equiangular lines in $\R^n$. Denote this quantity by $M(n)$. The purpose of this paper is to prove that there are no 76 equiangular lines in $\mathbb{R}^{19}$. Since 1973, $M(19)$ has been known to be between 72 and 76. The past four decades have seen no improvement on this bound. From Witt design, we can construct 72 equiangular lines in $\R^{19}$\cite[p. 148]{tay72}. The upper bound 76 is the relative bound for equiangular lines in \cite{lem73} and also a semidefinite programming bound in \cite{barg14}. After our results, $M(19)$ will be reduced to the range 72-75 and we conjecture that 72 is the maximum in $\mathbb{R}^{19}$.

The problem of determining $M(n)$ has been studied for almost seven decades, yet we still know very little.
Until recently the maximum number of equiangular lines in $\R
^n$ was known only for 35 values of the dimension $n$. We know $M(n)$ for most values of $n$ if $n \leq 43$. However, the cases for $n=14, 16, 17, 18, 19, 20$ and $42$ are still open. The history of this problem started with Hanntjes \cite{han48}  who found M(n) for $n = 2$ and $3$ in $1948$. Van Lint and Seidel \cite{lin66} found the
largest number of equiangular lines for $4 \leq n \leq 7$. In 1973, Lemmens and Seidel
\cite{lem73} determine M(n) for most values for $8 \leq n \leq  23$. Barg and Yu \cite{barg14} determine $M(n)$ for $ 24 \leq n \leq 41$ and $ n=43$. For other works on the bounds for equiangular lines, please see \cite{grea14}, \cite{bur15} and \cite{oy14}. 

We sketch below the approach of our main results. First, we will prove that if there exist 76 equiangular lines, the common angle has to be $\frac 1 3$ or $ \frac 1 5$. Then, by Lemmen-Seidel's bounds for equiangular lines with the angle $\frac 1 3 $, we know that $\frac 1 5$ is the only possible angle. Furthermore, if there exist 76 equiangular lines in $\R^{19}$ with the common angle $\frac 1 5$, it gives rise to an equiangular tight frame (ETF), which implies the existence of the strongly regular graph (75,32,10,16) \cite{wal09}. Azarija and Marc \cite{aza15} proved the nonexistence of the strongly regular graph (75,32,10,16). Therefore, there are no 76 equiangular lines in $\R^{19}$. 

By the classical treatment of strongly regular graphs (SRGs), the projection of the vertex set onto a
non-trivial eigenspace is a spherical 2-distance set and a 2-design \cite{car01}. The projection of srg(76,35,18,14) and its complement srg(76,40,18,24) will both form a spherical 2-distance set with inner product values $\pm \frac 1 5$ in $\R^{19}$. Namely, it gives rise to a 76 equiangular line set in $\R^{19}$ with the angle $\frac 1 5$ which contradicts  our main result. As a corrolary, we can show the nonexistence of these two SRGs. If the srg(76,30,8,14) or its complement srg(76,45,28,24) exists, there exist 76 equiangular lines in $\R^{57}$ with the  common angle $\frac 1 {15}$. By the existence of complementary ETFs, we will have 76 equiangular lines in $\R^{19}$ with the common angle $\frac 1 5$.  Therefore, these two SRGs do not exist.  Similar discussion also can prove that there is no 96 equiangular lines in $\mathbb{R}^{20}$ and the nonexistence of srg $(96,45,24,18)$ and srg$(96,38,10,18)$. In the last section, we discuss the connection between ETFs, SRGs, spherical few-distance sets and spherical designs.

\section{Prelimanaries}
A set of lines in $\R^n$ is called equiangular if the angle between each
pair of lines is the same. If we have $M$ equiangular lines in $\R^n$, then we will have  a set of unit vectors $\{x_i\}_{i=1}^M$ such that $|\langle x_i, x_j\rangle | =c$ for all $ 1\leq i \neq j \leq M$ , where $c$ is a positive constant. We call $c$ the \emph{common angle} of the equiangular lines. When the number of equiangular lines is large enough, the common angle will be the reciprocal of an odd integer. Neumann proved the following theorem.
\begin{thm}\label{thm:Neu}
 (Neumann \cite{lem73}) If we have $M$ equiangular lines in $\R^n$ and $M>2n$, then the common angle will be $\frac 1{2k-1}$, where $k \in \mathbb{N}$. 
\end{thm}  

Then, Larman, Rogers and Seidel proved a similar result for spherical two-distance sets.
A set of unit vectors $S = \{x_1, x_2, . . . \} \subset \R^n$ is called a spherical two-distance set if  $\langle x_i
, x_j \rangle \in
\{a, b\}$ for some $a, b$ and all $i \neq j$. The study of upper bounds of spherical two-distance sets can be found in \cite{barg13}. 
\begin{thm}\label{thm:LRS}
 (Larman, Rogers, and Seidel \cite{lar77}). Let $S$ be a spherical two-distance
set in $\R^n$. If $|S| > 2n + 3$ and $a > b$, then $b = \frac{ka-1} {k-1}$
for some integer $k$
such that $2 \leq k \leq \frac {1 + \sqrt{2n}}2$.
\end{thm}
The condition $|S| > 2n+ 3$ was improved to $|S| > 2n+ 1$ by Neumaier \cite{neu81}. 
If the spherical two-distance set gives rise to equiangular lines, then $a = -b$. So
Theorem \ref{thm:LRS} implies that $a = \frac 1 {2k-1}$, where $k \in \mathbb{N}$ , which is the statement of the Neumann
theorem in \cite{lem73}. The assumption of Theorem \ref{thm:LRS} is more restrictive than that of Neumann's theorem, but in return we obtain an upper bound on k. For instance when $n=19$, the common angle has to be  $\frac 1 3$ or  $\frac 1 5$ which cannot be deduced from Neumann's Theorem.

A finite collection of vectors $S = \{x_i \}_{i=1}^M \subset  \R^n$
is called a finite frame for the Euclidean space
$\R^n$ if there are constants $0 < A  \leq B < \infty$ such that for all $x \in \R^n$
$$ A||x||^2  \leq \sum_ { i=1}^M  |\langle x, x_i \rangle |^2  \leq  B||x||^2.$$
If $A = B$, then $S$ is called a tight frame. Benedetto and Fickus \cite{ben03} introduced a useful parameter of the frame, called the {\em frame potential}. For our purposes it
suffices to define it as $  FP(S)=\sum_{i,j=1}^{M}|\langle
x_i,x_j\rangle|^2.$ We can derive the lower bounds of frame potential and the minimizers are tight frames.

\begin{thm}{\cite[Theorem.6.2]{ben03}}\label{thm:BF}
   If $S$ is a set of unit vectors $\{x_i\}_{i=1}^M$ in $\R^n$ and $M>n$, then
  \begin{equation}\label{eq:fp}
    FP(S)\ge \frac {M^2} n
  \end{equation}
  with equality if and only if $S$ is a tight frame.
\end{thm}

If the set $S$ is a tight frame and equiangular, i.e. $|\langle x_i, x_j \rangle |=c$ for all $ i \neq j$, then $S$ is called an equiangular tight frame (ETF). ETFs have many nice properties. For instnace, they are Grassmanian frames \cite{str03} and attain the classical \emph{Welch bound}\cite{wel74}. The \emph{Welch bound} is the famous lower bound on the coherence.

\begin{thm}
If we have  a set of unit vectors $S=\{x_i\}_{i=1}^M$ in $\R^n$, then 
$$ \max_{ i \neq j} |\langle x_i, x_j\rangle | \geq  \sqrt{\frac{M-n}{n(M-1)}},$$ 
where equality holds if and if $S$ is an ETF. 
\end{thm}

When $M$ and $n$ are given for an ETF, the common angle is determined as $\sqrt{\frac{M-n}{n(M-1)}}$. We use ETF(n, M, c) to denote $M$ points ETF in $\R^n$ with the common angle $c$.  
ETFs are closely related to strongly regular graphs (SRGs) which form the main source of their constructions. A regular graph of degree $k$ on $v$ vertices is called strongly regular if every two adjacent vertices
have $\lambda$ common neighbors and every two non-adjacent vertices have $\mu$ common neighbors.
Below we denote such strongly regular graph by srg$(v, k, \lambda, \mu)$ .

Waldron \cite{wal09} proved that the existence of ETFs is equivalent to the existence of SRGs with certain parameters. 
\begin{thm} \cite[Corollary 5.6]{wal09}
There exists an equiangular tight frame of $M>n+1$ vectors for $\R^n$ if and only there exists a strongly regular graph $G$ of the type
$$
\text{srg}(M-1,k,\frac{3k-M}2,\frac k 2 ),\quad  k=\frac 1 2 M- 1 + (1-\frac M {2n}) \sqrt{\frac{n(M-1)}{M-n}}.
$$
\end{thm}

Consequently, we have the following lemma.
\begin{lem} \label{lem:srg75}   
The srg$(75,32,10,16)$ exsits if and only if ETF$(19,76,1/5)$ exists. 
\end{lem}
Furthermore, we can have two other SRGs connected to the existence of ETF(19,76,1/5).
\begin{lem} \label{lem:srg76}
If either srg$(76,30,8,14)$ or srg$(76,35,18,14)$ exsits then ETF$(19,76,1/5)$ exists. 
\end{lem}

Lemma \ref{lem:srg76} is a new result and we have two different ways to prove it. 
The first approach is based on the fact that the projection of the vertex set of an SRG onto a
non-trivial eigenspace is a spherical 2-distance set and a spherical 2-design \cite{car01}. Every spherical 2-design is a tight frame. Therefore, the projections of SRGs are \emph{two-distance tight frames} which have been discussed in \cite{bgoy15}. We define the notion of general spherical $t$-designs as follows. 

\begin{defn}\cite{del77b}
Let $\text{Harm}_t(\R^n)$ be the set of homogeneous harmonic polynomials of degree $t$ in $\R^n$. Let $t$ be a natural number. A finite subset $X$ of the unit sphere $S^{n-1}$
is called a spherical t-design if 
$$\sum _{x \in X} f(x) = 0,  \quad \forall f(x) \in \text{Harm}_j (\R^n), 1 \leq j  \leq t.$$ 
\end{defn}
We are interested
in the minimum cardinality of a spherical design when $t$ and $n$
are given. Delsarte, Goethals and Seidel \cite{del77b} proved
that the cardinality of a spherical $t$-design $X$ is bounded
below,
$$
|X| \geq   \binom {n+e-1}{n-1} + \binom {n+e-2}{n-1}, \quad |X| \geq 2  \binom {n+e-1}{n-1}
$$
for $t=2e$ and $t=2e+1$, where $e \in \mathbb{N}$. The
spherical $t$-design is called tight if the above bounds are
attained. If $X$ is a tight spherical $2s$-design, it is
immediately a maximum spherical $s$-distance set attaining
the linear programming bound in \cite{del77b}.
Also, $X$ is a spherical 2-design if and only if X is a tight frame with the center of mass at
the origin \cite{car01} \cite[Chapter 1]{yu14}. 

\begin{thm}\cite{del77b}\cite[Theorem 5.3 and  Proposition 5.1 ]{car01} \label{thm:2d2s}
Let G be a connected strongly regular graph which is not complete
multipartite, and let X be the projection of the vertex set of G onto a
non-trivial eigenspace, re-scaled to lie on the unit sphere. Then
X is a spherical two-distance set and a spherical 2-design.
\end{thm}

Notice that an SRG has two non-trivial eigenspaces. Therefore, every SRG gives rise to two different spherical two-distance sets which are also spherical 2-designs. 
\begin{example}
 If srg$(76,35,18,14)$ exists, then we will have two different spherical two-disance sets and 2-designs. The first, $S_1$, has 76 points in $\R^{19}$ with inner product values $\pm \frac 1 5$ and the second, $S_2$, has 76 points in $\R^{56}$ with inner product values $\frac {-3}{35}$ and $\frac 1 {20}$. Both of them are also spherical 2-designs.
\end{example}
$S_1$ gives rise to 76 equiangular lines in $\R^{19}$. Therefore, if srg$(76,35,18,14)$ exists, then ETF$(19,76,1/5)$ exists. By \cite[Proposition 3.1]{bgoy15}, the shifted 2-design of $S_2$ is a two-distance tight frame with inner product values $\pm \frac 1 {15}$ in one higher dimension ($\R^{57}$), i.e. it is an ETF(57,76,1/15). By the results of complementary equiangular tight frame \cite[Corollary 3.2]{cas13}, an ETF with $M$ elements in $\R^n$ exists if and only an ETF with $M$ elements in $\R^{M-n}$ exists. Therefore, ETF(57,76,1/15) exists if and only if ETF(19,76,1/5) exists.

Applying Theorem \ref{thm:2d2s} again, the projection of srg(76,30,8,14) results in an ETF(57,76,1/15) and a spherical two-distance set with 76 points in $\R^{18}$ with inner product values $\frac {-4}{15}$ and $\frac{7}{45}$. The shifted 2-design \cite{bgoy15} of the latter case has inner product values $\pm \frac 1 5$ in $\R^{19}$, i.e. it gives rise to ETF(19,76,1/5). Both of the projections of srg(76,30,18,14) give rise to ETF(19,76,1/5). 

To summerize, if either srg(76,30,8,14) or srg (76,35,18,14) exsits, then ETF(19,76,1/5) exists. Lemma \ref{lem:srg76} follows.

The second approach is as follows : if the Gram matrix of an ETF has the regular property (i.e. each row has the same number of $c$), we can use two different SRGs to construct the Gram matrix of an ETF. Conversly, if an ETF has the regular property, based on the tight frame conditon for Gram matrix, it gives rise to two different SRGs. This approach has been discussed in \cite{yu14}, \cite{bgoy15} and \cite{fjg15}. Following this approach, we can use the adjacency matrix of  srg(76,30,8,14) or srg (76,35,18,14) to construct the Gram matrix of ETF(19,76,1/5). We work indepedently and notice that \cite[Corollary 4.5]{fjg15} implies the same result for Lemma \ref{lem:srg76}.

\begin{thm} \cite[Corollary 4.5]{fjg15}
 If there exists an srg$(v, k, \lambda, \mu)$ with $v = 4k - 2\lambda - 2\mu$ then there exists an
srg$(v - 1, k \frac{v-2k}{v-2k-1}, \frac{3k-v}2 +\frac{3k}{2(v-2k-1)} , \frac k 2 \frac{v-2k} {v-2k-1}).$
\end{thm}

\section{Main results}\label{sec:main}
In general, there are no constraints on the common angle of equiangular lines. However, if there are 76 equiangular lines in $\R^{19}$, the common angle has to be $\frac 1 5$. 
\begin{lem} \label{lem1}
If there are 76 equiangular lines in $\R^{19}$, then the common angle of those equiangular lines is $\frac 1 5$.
\end{lem}
\proof
By Theorem \ref{thm:LRS}, since $76 > 2 \cdot 19+3$, then the common angle of 76 equiangular lines in $\R^{19}$ has to be  $\frac 1 3$ or  $\frac 1 5$. By Theorem 4.5 in \cite{lem73}, we know that if the common angle of equiangular lines is $\frac 1 3$ and $n \geq 15$, then the upper bound for such equiangular lines is $2n-2$. Since $76 > 2 \cdot 19-2$, then the common angle cannot be $\frac 13$.   \qed

\begin{lem} \label{lem2}
If there are 76 equiangular lines in $\R^{19}$ with the common angle $\frac 1 5$, then ETF(19,76,1/5) exists.
\end{lem}

\proof
If there are 76 equiangular lines in $\R^{19}$ with the common angle $\frac 1 5$, then there exists a set of unit vectors $S= \{x_i\}_{i=1}^{76} \in \R^{19}$ such that $|\langle x_i,x_j\rangle| = \frac 1 5$ for all $ 1 \leq i \neq j \leq 76$. Then, 
$$ FP(S) =\sum_{i,j=1}^{76}|\langle
x_i,x_j\rangle|^2 = 76+ 76 \cdot 75 \cdot (\frac{1}5)^2  = 76 \cdot 4 = \frac {76^2}{19}.$$
By  Theorem \ref{thm:BF}, since FP($S$) attains equality  \eqref{eq:fp}, then $S$ is an tight frame. Namely, $S$ is an equiangular tight frame, and hence ETF(19,76,1/5) exists. \qed 

\begin{thm}\label{main}
There are no 76 equiangular lines in $\R^{19}$. 
\end{thm}
\proof
By Lemma \ref{lem1} and \ref{lem2}, if there are 76 equiangular lines in $\R^{19}$, then ETF(19,76,1/5) exists. Furthermore, by Lemma \ref{lem:srg75}, there exists srg(75,32,10,16). However, this contradicts Azarija and Marc's result that there is no srg(75,32,10,16) \cite{aza15}.  
\qed

\begin{cor}
srg$(76,30,8,14)$ and srg$(76,35,18,14)$ do not exist. 
\end{cor}
\proof 
By Lemma \ref{lem:srg76}, if either srg(76,30,8,14) or srg(76,35,18,14) exists, then ETF(19,76,1/5) exists, i.e. there are 76 equiangular lines in $\R^{19}$. It contradicts Theorem \ref{main}.    \qed 

A. V. Bondarenko, A. Prymak and D. Radchenko proved the nonexistence of srg(76,30,8,14) \cite{bon14}. Here, we use the connetion between an SRG and a two-distance set to offer an alternative proof. For srg(76,35,18,14), the paper \cite{aza15} indicated that the proof is obtained from personal communication with  Haemers. Here, we offer the prove by the notion that the sphere embedding of an SRG will obtain a sphericla two-distance set which is also a spherical 2-design. 

\begin{rem}\label{rem}
Recently, nonexistence of the srg $(95, 40, 12, 20)$ is proved in \cite{aza16}. Similar discussion also can prove that there are no 96 equiangular lines in $\mathbb{R}^{20}$ and the nonexistence of srg $(96,45,24,18)$ and srg$(96,38,10,18)$.  
\end{rem}  
     
\section{Discussion}
We are interested in determining the maximum number of equiangular lines in $\R^n$. For $n=2,3,5,6,7,21,22,23$ and $43$, the maximum equiangular lines are ETFs. Previously, $n=19$ and 20 are conjectured the existence of ETFs to attain the upper bounds 76 and 96 respectively. However, we prove the case $n=19$ not attained. For n=20, there are three SRGs connected to the existence of ETF(20,96,1/5).  Using the same ideas in Lemma \ref{lem:srg75} and \ref{lem:srg76}, we can show that  srg(95,40,12,20) exists if and only if ETF(20,96,1/5) exists and if srg(96,45,24,18) or srg(96,38,10,18) exists, then ETF(20,96,1/5) exists. However, by \cite{aza16}, we know none of them exist. (Ref: Remark \ref{rem}) 

We note that for several of the sets of parameters that correspond to open
cases in Table \ref{table:ETFs}, their cardinality matches the
best known upper bound on the size of equiangular line set in that
dimension (the semidefinite programming, or SDP, bound of
\cite{barg14}). If we know the existence of any SRGs in Table \ref{table:ETFs}, we will obtain new results for maximum equiangular lines in that dimension. Specifically, this applies to $n=42,45,46.$
For instance, in the case of $n=42$ the SDP bound gives $M=288$
and $c=1/7$ (it is not known whether a set of 288 equiangular
lines in $\R^{42}$ exists). Using our approach, we observe that
such a set could be constructed from srg(287,126,45,63), srg(288,140,76,60) and
srg(288,164,100,84). Unfortunately, neither of them is
known to exist (or not). Notice that in table \ref{table:ETFs}, we know the nonexistence of srg(540, 308,190,156).  However, this is not sufficient to show the nonexistence of ETF(45, 540, 1/7). Therefore, the existence of 540 equiangular lines in $\mathbb{R}^{45}$ remains an open question. Furthermore, we want to connect the notion of tight spherical 5-designs with ETFs. 

\begin{thm}[Gerzon]\cite{lem73}
If there are $M$ equiangular lines in $\mathbb{R}^n$, then
  \begin{equation}\label{eq:gerzon}
M \leq \frac{n(n+1)}{2}
  \end{equation}
\end{thm}

 Currently, Gerzon's bounds are attained only for
$n=2,3,7,$ and $23.$
Note that if there are equiangular lines attaining the Gerzon bound, then the common angle is $\frac1{\sqrt{(n+2)}}$ 
\cite[Thm.3.5]{lem73}. 

\begin{thm}\cite[Theorem 5.12]{del77b}
If $S$ is a tight spherical 5-design in $\R^n$, then $|S|=n(n+1)$ and the inner product values between each pair of points in $S$ are -1 and the zeros of the polynomial $C_2(x)=1+ \frac{ (n+2)(nx^2-1)}2.$
\end{thm}
The zeros of $C_2(x)$ are $\pm \frac{1}{\sqrt{n-2}}$. If $S$ is  a spherical tight 5-design, then $S$ is antipodal and inner product values are -1 and $ \pm \frac{1}{\sqrt{n+2}}$. Therefore, tight spherical 5-designs will give arise to $\frac{n(n+1)}2$ equiangular lines in $\R^n$ and vice versa. By Neumann's Theorem, when $n>2$, $\frac 1{\sqrt{n+2}}=  \frac 1 {2m+1}$, where $m \in \mathbb{N}$. Therefore, $n$ has to be an odd square minus two, i.e. $n= (2m+1)^2 -2$ for some positive integer $m$. The existence of a tight spherical 5-design in $\R^n$ is equivalent to the existence of an ETF$(n,\frac{n(n+1)}2, \frac{1}{\sqrt{n+2}})$. Futhermore, such ETFs minimize potential energy for each completely monotonic potential function, i.e. they are universal optimal codes \cite{ck07}. For instance, for the cases of $m=1$ and $2$, ETF(7,28,1/3) and ETF(23,276,1/5) form very nice configurations in the corresponding dimension. Based on the results of E. Bannai, A. Munemasa, and B. Venkov \cite{ban04} and Nebe and Venkov \cite{neb12}, there are no tight 5-designs in $\R^n$, where $n=(2m+1)^2-2$ with an infinite sequence of values of $m$ that begins with $m =3, 4, 6, 10, 12, 22, 38, 30, 34, 42, 46.$ For $m=3$ and $4$, we that ETF(47,1128,1/7) and ETF(79,3160,1/9) do not exist. Using the same ideas in Lemma \ref{lem:srg75} and \ref{lem:srg76}, we can show the nonexistence of srg(1127,640,396,320), srg(1128,644,400,324), srg(1128,560,316,240), srg(3159,1408,1064,702),  srg(3160,1575,870,700) and srg(3160,1755,1050,880). Note that the first two cases are known to not exist in Brouwer's table \cite{bro15}. 

By \cite[Theorem 5.11]{del77b}, tight spherical 4-designs in $\R^n$ are the maximum spherical two-distance sets with inner product values $\frac{-1 \pm\sqrt{n+3}}{n+2}$. Also, $n$ has to be odd square minus three. For instance, if $n=22$, the tight spherical 4-designs in $\R^{22}$ are the maximum spherical two-distance set in $\R^{22}$ with 275 points and inner product values are $\frac 1 6$ and $-\frac 14$. Such a spherical two-distance set can be obtained from the projection of srg(275,112,30,56) and the projection of 276 equiangular lines in $\R^{23}$. This observation may offer another point of view in recognizing that tight spherical 5-designs in $\R^n$ are equivalent to tight spherical 4-designs in $\R^{n-1}$. There are also more connections between tight spherical designs of hamonic index $T$ and spherical few-distance sets in \cite{ban09}, \cite{ban13}, \cite{oy14}, and \cite{zhu15}. In this discussion, we like to relate different notions for mathematicians who are interested in frame theory, SRGs, equiangular lines, spherical few-distance sets, spherical t-designs and some related topics in algebraic combinatorics.

\begin{table}\begin{center}
{\small \begin{tabular}{|c|c|c|c|c|} \hline
  $n$ & $N$ & $c$ & comments\\
\hline

42&288&1/7& srg(287,126,45,63)(o)\\ &&&  srg(288,140,76,60) (o)\\&&& srg(288,164,100,84) (o)\\
45&540&1/7 &  srg(539,234,81,117)(o)\\ &&& srg(540,266,148,114) (o)\\&&& srg(540,308,190,156) (N)\\
46&736&1/7 & srg(735,318,109,159)(o)\\ &&& srg(736,364,204,156) (o)\\ &&& srg(736,420,260,212) (o)\\
\hline
\end{tabular}}
\vspace*{.1in}\caption{Parameter sets of possible maximum ETFs. The
label `o' means that the existence of an SRG with these parameters
is an open problem. `N' means that the srg does not exist.}\label{table:ETFs}
\end{center}
\end{table}

\section*{Acknowledgements.}
The author would like to thank Alexander Barg and Alexey Glazyrin whose suggestions and comments were of inestimable value for this paper. The author also thanks Ye-Kai Wang and Aditya Viswanathan for their valuable comments.

\providecommand{\bysame}{\leavevmode\hbox
to3em{\hrulefill}\thinspace} \providecommand{\href}[2]{#2}
\bibliographystyle{amsalpha}

\end{document}